\documentclass[11pt, oneside]{article}   	
\usepackage{geometry}                		
\usepackage{graphicx}	
\usepackage{hyperref}			

\usepackage{amsmath}
\usepackage{amsfonts}
\usepackage{amssymb}
\usepackage[version=4]{mhchem}
\usepackage{stmaryrd}

\usepackage{epstopdf}
\usepackage[utf8]{inputenc}

\title{Pierre Cartier: A Visionary Mathematician}
\author{Alain Connes\footnote{Coll\`ege de France} \, \& Joseph Kouneiher\footnote{C\^{o}te d'Azur University}}

\begin{document}
\maketitle

\abstract{This article reflects on the life and mathematical contributions of Pierre Cartier, a distinguished figure in 20th- and 21st-century mathematics. As a key member of the Bourbaki collective, Cartier played a pivotal role in the formalization and modernization of mathematics. His work spanned fields such as algebraic geometry, representation theory, mathematical physics, and category theory, leaving an indelible mark on the discipline. Beyond his technical achievements, Cartier was celebrated for his intellectual generosity and humanistic approach to science, shaping not only mathematical thought but also the broader cultural understanding of the field.
}

\section*{Introduction}

 Pierre Cartier passed away on August 17th, 2024, leaving behind an intellectually rich and profoundly human legacy. The first word that comes to mind when thinking of Pierre is universality. He was an extraordinary mathematician, whose remarkable intuition impressed many. Alexandre Grothendieck, himself a legend in the field, recognized in R\'ecoltes et Semailles Cartier's unparalleled intuition, capable of penetrating the most varied subjects with remarkable clarity and depth.

 Pierre dedicated his life to mathematics, to which he contributed a vast treasure, both through his personal discoveries and the ideas he generously shared with the community.

 Pierre Cartier was born on June 10, 1932, in Sedan, France. After completing his secondary education in Sedan and preparatory classes at Lyc\'ee Saint-Louis in Paris, he was admitted to the \'Ecole Normale Sup\'erieure in 1950 through the mathematics competition, and defended his thesis in 1958 under the supervision of Henri Cartan. Officially, at the beginning, his thesis advisor was Roger Godement. However, he felt more inspired by the work of Cartan and especially Andr\'e Weil, so he changed his research topic ("Derivations and Divisors in Algebraic Geometry").

\begin{figure}
\begin{center}
 \includegraphics[width=12cm]{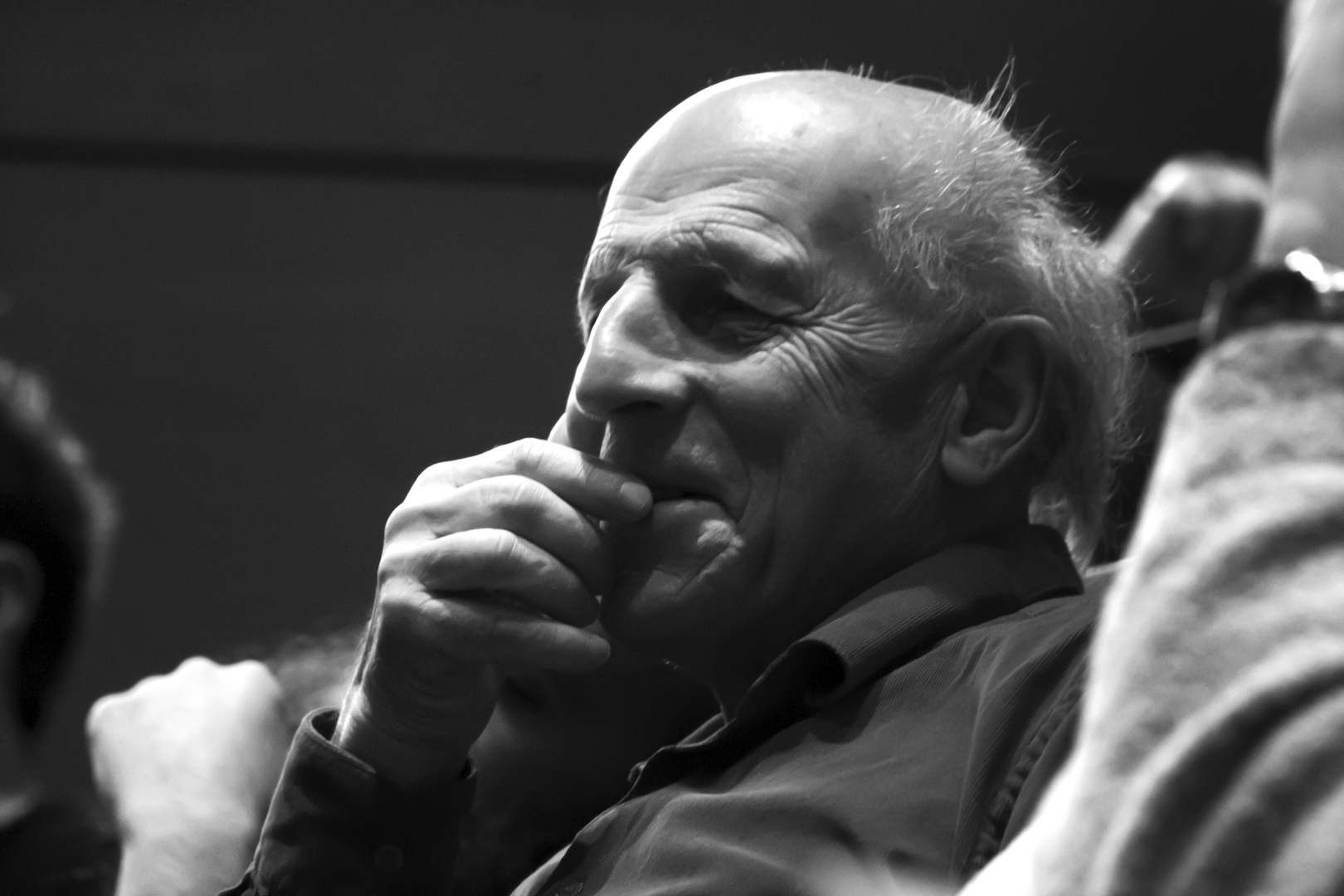}

 \label{fig:cartier-pic}
 \caption{Pierre Cartier}
 \end{center}
\end{figure}

\begin{quote} The best result of my thesis, the duality of abelian varieties \cite{cartier1}, was a problem posed by Andr\'e Weil in his book on abelian varieties and algebraic curves.
\end{quote}

 He became a member of the Bourbaki group in 1955, at the age of 23, and remained so until 1983. After holding a professorship in Strasbourg during the 1960s, he was appointed Research Director at the CNRS. He was a permanent visitor at the Institut des Hautes \'Etudes Scientifiques (IHES) in Bures-sur-Yvette, where he left a lasting impact, and served as a visiting professor or researcher at Princeton, the \'Ecole Polytechnique, and various other institutions. From 1988 to 2002, he was a professor at the \'Ecole Normale Sup\'erieure. He is best known for his work in algebraic geometry, representation theory, homological algebra, and category theory. His interest in mathematics emerged very early, particularly after reading Hermann Weyl's work on group theory and quantum mechanics at the age of 18:

\begin{quote} 
It was undoubtedly my reading, at the age of eighteen, of Hermann Weyl's classic work on group theory and quantum mechanics that had the greatest influence on the course of my subsequent mathematical research. In it, I discovered how the beauty of a fully developed mathematical theory can be combined with a profound physical necessity, and the importance of the connections drawn between seemingly unrelated fields.

Through contact with eminent scholars who taught mathematics at the \'Ecole Normale Sup\'erieure, and later with Nicolas Bourbaki, I definitively oriented myself towards 'pure' mathematics, without disregarding the value of analogies drawn from the physical sciences, the importance of applications, or the assistance provided by a well-controlled intuition.

If the list of my publications gives the impression of incessant wandering--from algebraic geometry to probability theory, from algebra to quantum mechanics, from group theory to numerical analysis, or from number theory to combinatorics--it is because I have strived to perceive the multiple consequences of a few general ideas that I have constantly revisited. I have been more interested in the discovery of new truths than in the meticulous cultivation of a field already well explored (\cite{notices}, p. 1).
\end{quote}

 For years, Pierre was a pillar of the Bourbaki collective, playing a key role in transforming 20th-century mathematics and helping to propagate the language and methods that have become standard in modern mathematical research. He was one of the rare individuals who mastered the full breadth of Bourbaki's work, reflecting his deep belief in the power of structure and generality in mathematical thinking. Cartier dedicated a significant portion (one-third, according to him) of his mathematical activity to drafting Bourbaki's books. For instance, chapters 4 to 6 of the famous books on Lie groups and Lie algebras benefited from Cartier's idea to base everything related to root systems, Weyl groups, etc., on symmetries with respect to hyperplanes in a vector space. He also delivered about 40 Bourbaki seminars, reflecting his commitment to the collective and to the transmission of knowledge.

 Pierre Cartier's contributions are as vast as they are varied, spanning fields such as algebraic geometry, number theory, category theory, distribution theory, and mathematical physics. Through the concepts he created, many of which bear his name, such as Cartier divisors in algebraic geometry, he provided the correct conceptual framework for understanding key notions. For instance, in the case of divisors, his idea--derived from Henri Cartan's work in analytic geometry--of a local sheaf-theoretic formulation applies even in situations\footnote{For instance in tropical geometry}  that would be inaccessible to a more naive notion of divisors. He was the first to uncover the special properties of the de Rham complex in finite characteristic, offering groundbreaking insights into its behavior in this setting, and introduced the Cartier operation in the differential calculus of algebraic varieties in positive characteristic. This revealed a radically different scenario compared to characteristic zero.

 Cartier's work in the theory of commutative formal groups represents a foundational contribution to algebraic geometry and number theory. His introduction of \emph{Cartier's formal groups} provided a deep and flexible framework for studying algebraic curves, group schemes, and their associated deformation theory. In essence, Cartier formal groups offer a language for understanding the infinitesimal behavior of varieties and schemes over fields, rings, and more general bases.

 Regarding his involvement in Grothendieck's work, Cartier wrote:

\begin{quote}
All these concepts such as {\it Cartier divisors}, the {\it Cartier operation on differential forms}\footnote{This operation allowed Ogus to resolve Katz's conjecture in the general case, predicting certain inequalities on p-adic valuations of Frobenius. The construction, now known as the Cartier isomorphism, is a reformulation by Katz of Cartier's work published in a 1957 note in Comptes Rendus. This isomorphism, which relates the components and cohomology groups of the de Rham complex of a smooth variety in positive characteristic, would come to dominate all of differential calculus in characteristic $p$ or mixed characteristic for years to come.}, {\it the Cartier duality} of group schemes with finite group schemes, and particularly {\it inseparable descent }translate seamlessly into \newline Grothendieck's new theory of schemes. By combining his methods with mine, Grothendieck obtained a very simple proof of the completeness theorem for linear systems of divisors (\cite{notices}, p.9).
\end{quote}

\noindent and he adds :

\begin{quote}
In 1966, I glimpsed the possibility of extending Dieudonn\'e's structure theorems to this new framework. My aim was initially to generalize the construction of one-dimensional formal groups introduced by Lubin and Tate, and to express it independently of the choice of a coordinate. I first demonstrated that the essential property was the existence of a lifting of the Frobenius operator, and I was able to obtain a complete description in terms of modules for the corresponding class of formal groups. Using these results, I was able to fully resolve a series of problems formulated by Grothendieck\footnote{In his article, On Some Points of Homological Algebra (Toh\^{o}ku), after stating the lack of a satisfactory theory for multiplicative structures in Homological Algebra that meets the necessary level of generality and simplicity, Grothendieck adds a footnote acknowledging Pierre Cartier's contribution: ``Mr. P. Cartier has just found a generally satisfactory formulation for multiplicative structures in Homological Algebra".} in connection with his theory of crystals.
\end{quote}

\section{Hopf algebras, group theory}
\noindent  Cartier had a guiding thread in his work that revolved around groups and Hopf algebras \cite{cartier12, cartier010}. He viewed Hopf algebras as a powerful tool for understanding various algebraic phenomena. His research often focused on their applications in algebraic topology, representation theory, and mathematical physics. The duality in Hopf algebras, where every algebraic structure has a corresponding co-algebra structure, resonated with Cartier's broader mathematical philosophy, which emphasized the interconnectedness of different fields of mathematics.
\begin{quote}
I have never abandoned group theory: for me, it remains the central point of everything I have done. In my opinion, the book that everyone should have read is Group Theory and Quantum Mechanics by Hermann Weyl. It's a text that I still read today with the same interest. Group theory has the advantage of allowing one to do physics in the way Weyl proposed, so it's an important tool in physics. It's also an important tool in geometry, of course, after the works of \'Elie Cartan (the father of Henri Cartan), and even before that. It's a crucial tool in arithmetic, as demonstrated by Andr\'e Weil. For me, it's as if I have a fortress, the fortress of groups: I can explore here and there, enter through another door, but in the end, I always return to my fortress\footnote{personal communication}.
\end{quote}

 His early interests and original results focused on characterizing the enveloping algebra of Lie algebras through the properties of its filtration and coproduct, establishing the algebraic foundation for the local theory of Lie groups and its extension to formal groups. Another central idea is his geometric theory of Weyl groups. His perspective is that the classical theory should be understood as the study of groups generated by symmetries relative to linear varieties. At a Bourbaki congress in Pelvoux\footnote{The Bourbaki congresses are meetings of the members of the Bourbaki group in which participants discuss writings by one of them on selected topics which after many versions will become the Bourbaki books}, Pierre Cartier introduced the axioms of root systems, which play a key role in the chapters on Lie groups and Lie algebras.

\noindent One of the recurring themes in his research was the relationship between the representations of a Lie group and those of its Lie algebra, specifically the problem of infinite-dimensional unitary representations. Along with Dixmier  \cite{dixmier}, he extended Harish-Chandra's results to the case of arbitrary Lie groups in a Banach space. These results caught the attention of mathematical physicists who, following in the footsteps of Hermann Weyl and von Neumann, were analyzing commutation relations in quantum mechanics. Cartier introduced distribution vectors for infinite-dimensional representations, which were subsequently used in functional analysis and quantum field theory. In fact, Ed Nelson was able, in 1959, to remove an unnecessary hypothesis in  \cite{dixmier}  by using entirely new methods in the theory of processes: stopping times and the strong Markov property. This triggered Pierre's interest in probability theory. In 1964, he reported on the early work of the Russian school (Gelfand, Milnos) on stochastic processes in infinite-dimensional vector spaces. He introduced the class of "standard spaces" into the theory, in which measure theory does not present the "monstrosities" that had so hindered the development of process theory. His student, Fernique, proved in his thesis that "all" functional spaces are standard. Cartier observed that these standard spaces were none other than non-metrizable Lusin spaces and that the entire classical theory of Lusin and Souslin never relied on the metrizability hypothesis, as initially assumed. His collaboration during this period with P. A. Meyer and L. Schwartz culminated in the Bourbaki volume that concluded the series on integration and measure theory.

\noindent Cartier then turned his attention to the Markovian properties of random functions (or distributions) with multiple parameters. His work resolved the controversy between Paul L\'evy and McKean, favoring the latter, and uncovered the flaw in Paul L\'evy's reasoning. These same random distributions found applications in the constructive theory of quantum fields, following the work of Nelson and Cartier himself. 
In his work on the path integral formulation of quantum mechanics and quantum field theory, Cartier, in collaboration with Cecile Dewitt, sought to place this crucial concept on solid mathematical ground. By integrating methods from algebra, geometry, and analysis, Cartier helped to clarify the path integral's use and its connections with other areas of mathematics, such as stochastic processes and functional integration.

\noindent We shall return in Section 1.3 to another major contribution of Pierre to mathematical physics, namely his cosmic Galois group in renormalization theory.

\subsection{On the theory of Formal groups }

 \noindent Among Cartier's original contributions to group theory are his results related to formal groups, following Dieudonn\'e's publications in 1954 and 1959, which established the theory of formal groups. This work marked a return to the roots of Lie theory by studying power series that express multiplication near the origin of a Lie group. Dieudonn\'e, especially in the case of commutative groups, obtained a series of new results.

\noindent Cartier's aim, building on H. Cartan's results in algebraic topology, was to reformulate Dieudonn\'e's theory in terms of linear duality. Cartier emphasized the importance of the coproduct in Dieudonn\'e's "hyperalgebra," providing an intrinsic formulation of the notion of a formal group, establishing its equivalence with that of filtered bialgebras, and extending the framework to obtain the kernels of homomorphisms.

\noindent Building on Dieudonn\'e's work, Cartier observed that $p$-adic matrices were what Weil had lacked to complete the theory of abelian varieties.

 \begin{quote}   
    Upon reading the unpublished manuscripts on formal groups that Dieudonn\'e had shared with me, I had the sudden insight that his $p$-adic matrices were precisely what Weil had been missing to complete the theory of abelian varieties he had just developed. At the same time, Grothendieck was beginning to completely revolutionize algebraic geometry (see \cite{notices}, p.1).
    \end{quote}

\noindent  Under Jean-Pierre Serre's guidance, Cartier resolved all the issues left unresolved by Weil, including the biduality of abelian varieties, the absence of torsion in the divisor class group (also resolved by Barsotti), and the representation of homomorphisms of abelian varieties, all while employing the technique of studying isogenies, especially in the inseparable case.

\noindent  However, these works gained their full significance only after the introduction of group schemes by Grothendieck. During this period, Cartier introduced several key concepts: Cartier divisors, the Cartier operation on differential forms, Cartier duality of finite group schemes, and most notably, inseparable descent. All these concepts were seamlessly integrated into Grothendieck's new theory of schemes. By combining his methods with those of Cartier, Grothendieck provided the proof of the completeness theorem for linear systems of divisors.

\noindent  In 1966, Cartier recognized the possibility of extending Dieudonn\'e's structure theorems:

\begin{quote}
    The essential progress was the introduction of curve modules associated with formal groups and the notion of a typical curve. This novel concept retrospectively explained the arithmetic properties of the classical exponential series and, in particular, the existence of the Artin-Hasse exponential. It also offered a flexible interpretation of the Dieudonn\'e module and allowed me to refine known theorems on Witt groups. The universal typical curve in a Witt group plays a central role in these demonstrations\footnote{Personal communication} \cite{notices}.
    \end{quote}

\noindent  Cartier's theory of commutative formal groups and his work on the typification of curves through the introduction of what are now known as Cartier's formal groups were pivotal \cite{cartier5, cartier7}. These results laid the groundwork for understanding the structure of formal groups and their role in number theory and algebraic topology. His introduction of the global Witt construction, an advanced technique that unified various local constructions into a coherent global framework, remains a cornerstone in the study of algebraic structures.

  \subsection{Quantum groups}

\noindent Cartier was deeply interested in the theory of quantum groups, which are deformations of classical groups that can be studied within the framework of Hopf algebras. This work not only broadened the theoretical landscape of algebra but also had significant implications for mathematical physics, particularly in the study of symmetries and integrable systems. Drinfeld quantum groups, for example, are crucial for understanding how algebraic structures encapsulate quantum symmetries.

\noindent Cartier expressed his deep attachment to Hopf algebras when describing the most important scientific moment of his life: the presentation he gave at the International Congress of Mathematicians in 1986 in Berkeley and the new lines of research it opened:

\begin{quote}
    At the congress, I presented on behalf of Drinfeld, who was invited but couldn't attend due to the Soviet regime. On the first day, the Russian president of the International Mathematical Union handed me Drinfeld's English manuscript, unsure of what to do with it. There was also a text from Manin, and he asked if I could replace one of them. After a brief thought, I chose Drinfeld's. He then informed me that the presentation was scheduled for that afternoon. I took the manuscript and said, "I'll try." I explained the situation to Kaplansky, the American president of the organizing committee, who helped by locking me in a room with sandwiches and coffee. I had six hours to prepare. Although I was familiar with Hopf algebras, the material was new, and by the time I gave the talk, 400 people had gathered to listen to Drinfeld, which meant me. In the following days, I distributed as many copies of the text as I could, made using simple methods.

 At that time, I was nearing the end of my mathematical research program and was in a period of waiting. I realized that quantum groups were something entirely new that could be approached with techniques I had previously used (like Lie groups and differential geometry), which I knew well. There was also motivation from physics, even if not always apparent, which reoriented my interests for a good ten years. In fact, I still organize a seminar at the \'Ecole Polytechnique called "Quantum Groups and Poisson Geometry," which continues this work, twenty years later\footnote{Personal communication}.

\end{quote}

\subsection{The Cosmic Galois Groups}\label{cosmic}

 In the realm of mathematical physics, Pierre Cartier had a long-lasting interest in quantum field theory \cite{cartier6, cartier8, cartier11}. He proposed the idea \cite{cartier10} of a "cosmic Galois group" underlying the symmetries of the renormalization process, while the traditional renormalization group could be viewed as a one-parameter subgroup of this larger cosmic Galois group. On this idea, he writes:

\begin{quote}
    My notion of a "cosmic Galois group" arose from my reflections on mathematical physics. I drew particular inspiration from the work of \cite{connes}, who redefined the renormalization process--a technique for eliminating infinities in divergent integrals linked to Feynman diagrams. They introduced a new group structure to describe the intricate relationships in these calculations.

    At the same time, my mathematical research has focused on series and integrals involving special numbers, like powers of $\pi$ and values of the Riemann zeta function \cite{cartier4}. These relationships are governed by a symmetry group that resembles Grothendieck's motivic Galois group. I began to see a striking analogy between this group and the Connes-Kreimer group, suggesting that they might be two variations of the same group, influencing both the mathematical and physical aspects of the problem.

    The motivic Galois group deals with the automorphisms of certain transcendental numbers, which are similar to the constants appearing in Feynman diagram computations. This observation led me to interpret the group in \cite{connes} as a symmetry group governing the fundamental constants of physics. In the Standard Model, these constants are often adjusted empirically, without much mathematical explanation. I believe that this group could express new symmetries among these constants, which might have significant implications for cosmology. My ultimate dream is to unite the ideas in \cite{connes} with the motivic Galois group, even if, for now, it's still an ongoing research program\footnote{personal communication} \cite{cartier10}.
\end{quote}

It took a lot of work and time to obtain a concrete realization of Pierre Cartier's idea of the Cosmic Galois Group as well as its action on the coupling constants of renormalizable physical theories. The first step, following the introduction of the Hopf algebra of Feynman graphs, was the realization that the renormalization process is, in fact (in the dimensional regularization scheme), identical to the Birkhoff decomposition in pure mathematics, which occurs in the classification of vector bundles on the sphere (see \cite{CK}). Once this result was obtained, it again took time, using the Riemann-Hilbert correspondence, to realize (\cite{CM}) that there exists a universal Hopf algebra $\mathcal H$ behind all these computations\footnote{introducing the category of ``equisingular flat connections," which is shown to be a Tannakian category, meaning it is equivalent to the category of modules over a certain pro-algebraic group}. This is the graded dual $\mathcal H$ of the universal enveloping algebra of the free graded Lie algebra generated by elements of degree $n$ for each positive $n$. Let $\mathbb U$ be the affine group scheme associated with $\mathcal H$, and let $\mathcal H$ and $\mathbb U^*=\mathbb U \rtimes \mathbb G_m$ be its semidirect product by the action of the multiplicative group given by the grading. It follows (see \cite{CM}, Corollary 1.107) that this group acts on the coupling constants of any renormalizable physical theory, thus providing the best model for the Cosmic Galois Group. It gives a concrete realization of the symmetry group of physical theories and has a tantalizing similarity with the motivic Galois group of mixed Tate motives (see \cite{CM}, Corollary 1.111).

Curiously, one of Cartier's most cited results is the so-called Cartier theorem, often discussed in the context of algebraic groups over fields of characteristic zero. This theorem, a key result in algebraic geometry and group theory, asserts that {\it algebraic groups over fields of characteristic zero are smooth}. It plays an important role in the study of algebraic groups, as smoothness is a fundamental property with significant geometric and arithmetic implications.

In this context, smoothness means that the algebraic group is a smooth variety--its tangent space at every point is well-behaved, which roughly implies that there are no singularities in the variety. This result is foundational in the theory of algebraic groups because smoothness is a desirable property in geometric structures, ensuring that algebraic groups can be analyzed using standard differential-geometric tools (as they behave ``well" locally) and that they exhibit favorable properties in representation theory.

The fact that this theorem is only mentioned as a footnote in one of Cartier's papers \cite{cartier109} is due to its relatively technical nature and the fact that it was not the central focus of his work at the time. Although the proof of smoothness for algebraic groups in characteristic zero was not a primary goal in that paper, it was still an important result for the completeness of his theory, which is why it is briefly noted.

\section{His many key ideas offered to others} 
 
\noindent Pierre Cartier's influence extends well beyond the concepts associated with his name, as his intellectual generosity fostered significant collaborations with his contemporaries. Notably, he provided Andr\'e Weil with key ideas for proving fundamental results in number theory using the local compactness of the topological ring of adèles. Cartier's insights certainly influenced Weil in basing his book, {\it Basic Number Theory}, on the local compactness of the ring of adèles, which contains a global field as a discrete and cocompact subfield.

\noindent Cartier also advised Grothendieck to use prime ideals instead of maximal ideals when defining the spectrum of a ring in algebraic geometry, as prime ideals are the only ones compatible with morphisms. In classical algebraic geometry, points of an algebraic variety correspond to maximal ideals of its coordinate ring, while prime ideals correspond to irreducible closed subsets. The spectrum of a ring, $\text{Spec}(A)$, therefore includes standard points associated with maximal ideals and generic points linked to irreducible closed subsets\footnote{Grothendieck's theory of schemes builds on these concepts, bridging classical geometry with modern spatial concepts. A common justification for identifying the points of an affine scheme with prime ideals (instead of just maximal ideals) is that a ring homomorphism $\phi : A \longrightarrow B$ does not always induce a well-defined map from the set of maximal ideals of $B$ to $A$ ; however, the inverse image of a prime ideal under such a homomorphism $\phi : A \longrightarrow B$ remains a prime ideal in $A$.
}. The idea of viewing irreducible subsets of an algebraic variety as “points” can be traced back to Italian algebraic geometers in the early 20th century.

\vspace{0.2cm}

 Another clarification brought by Cartier to Grothendieck's work also concerns the concept of schemes. In algebraic geometry, schemes are mathematical objects that generalize algebraic varieties. Grothendieck introduced this concept in the late 1950s, defining schemes as entities that encode information about the solutions to polynomial equations through commutative rings. However, to address the complexity of these objects, Cartier proposed a more abstract and powerful perspective: he suggested viewing Grothendieck's schemes as functors from the category of commutative rings to the category of sets. This innovative perspective significantly impacted algebraic geometry by providing a more abstract and robust understanding of schemes.

The significance of this proposal lies in its universality and generalization. This functorial approach allows the concept of a scheme to be extended to more general situations, thereby facilitating the study of the properties of these objects in various contexts. This perspective also paves the way for applications in other areas of mathematics, including topology, number theory, and even theoretical physics, where similar constructions often prove useful. 
 It is worth noting that Chevalley had already used this idea in his 1956 paper on finite Lie groups, in which he defined an algebraic group as a functor from rings to groups.

\vspace{0.2cm}

The introduction of group varieties in transcendental number theory by S. Lang followed a conjecture of Cartier, who asked whether it would be possible to extend the Hermite-Lindemann theorem from the multiplicative group to a commutative algebraic group over the field of algebraic numbers. This is the result that Lang proved in 1962. At that time, there were a few transcendence results (by Siegel and Schneider) concerning elliptic functions and even Abelian functions. However, Lang's introduction of algebraic groups in this context marked the beginning of several important developments in the subject.

More precisely, in \cite{lang}, Lang proved Cartier's conjecture, which states that if $ G $ is an algebraic group over a number field $ K $ and $ \alpha \in \operatorname{Lie}(G)(K) $ is such that $ t \mapsto \exp _{G}(t \alpha) $ is not an algebraic function, then $ \exp(\alpha) $ is transcendental over $ K $. For $ G $ as a linear group, this reduces to the classical result concerning the exponential function. The novelty comes from the non-linear case; when $ G $ is an abelian variety, Lang's result represents a transcendence result for values of theta functions. Lang derived this theorem from his transcendence criterion, which generalizes the method of Gelfond and Schneider.

\section*{The other facets of Pierre Cartier}
 
 Beyond his technical contributions, Pierre Cartier was a passionate advocate for the philosophical and humanistic aspects of mathematics \cite{cartier13}. He recognized that mathematics is not merely a collection of theorems and proofs but a human endeavor that reflects creativity, beauty, and the search for truth. His writings and lectures often explored the broader implications of mathematical ideas, touching on their philosophical foundations and connections to other fields, such as physics and philosophy\footnote{His interest in the history and philosophy of mathematics is exemplified by the seminar he co-hosted for over 30 years.}.

Regarding his interest in the epistemology and philosophy of mathematics, he made the following remark:

\begin{quote} 
What attracts me to the epistemology of mathematics is understanding how mathematics is woven into civilization. I'm interested in how mathematical concepts emerge from societal concerns, reflecting the zeitgeist. One of my students specialized in the history of mathematics in China and worked on her thesis with me and my sinologist brother. Andr\'e Weil sparked my interest in the history of mathematics, teaching me to view past mathematicians -- like Euclid, Archimedes, Fermat, Euler, and Gauss -- as contemporaries. Hermann Weyl's writings also influenced me, prompting a desire to explore the relationship between physical and mathematical concepts from a philosophical perspective. Recently, I've focused on the history of categories, partly because I've contributed to their development and can draw on my memories. This area of mathematics closely aligns with philosophical issues, which I find more compelling than formal logic itself. The central question remains: what guarantees that mathematics conveys truth, and how does it do so consistently\footnote{personal communication}? \cite{fresan, cartier13}.
 \end{quote}

Meeting Pierre Cartier was an encounter with a personality who combined incredible intellectual rigor with an equally exceptional human generosity. His colleagues remember him as a man of great simplicity, whose altruistic pragmatism was as impressive as his vast mathematical knowledge. Anecdotes abound about his cycling exploits, demonstrating a physical vitality that seemed to defy time, and about  his talents in astronomy, a field in which he enthusiastically shared his knowledge.

On his approach to research, he says:

\begin{quote} 
Feynman once said that to be a genius, you just need to keep ten problems in mind and constantly look for solutions in everything around you. My approach is rooted in my innate curiosity, which was nurtured from my early studies in philosophy, physics, and mathematics. I always tackle multiple problems and methods simultaneously, allowing me to draw analogies between different subjects--a process that has often led to significant discoveries in my career.

My character was shaped by my father's boundless imagination and my grandmother's Alsatian common sense. This combination instilled in me an insatiable curiosity about people, travel, and reading on a wide range of topics. My wife, who had literary and musical interests, also opened up new perspectives for me, such as music, which I learned to appreciate thanks to her.

Today, I continue to cultivate this curiosity, for instance, by collaborating with musicians on projects inspired by Euler's classical works on music \cite{fresan, cartier13}. 
\end{quote}

Cartier's legacy also includes his role as a mentor and educator. He inspired countless students and colleagues with his enthusiasm for mathematics and his commitment to intellectual rigor. His generosity in sharing ideas, his encouragement of young mathematicians, and his ability to see connections across diverse areas of study have left a lasting impact on the mathematical community.

His influence extends far beyond his published work. He was known for his ability to ask deep, probing questions that often led to new avenues of research. His insights, whether shared in written form or in conversation, have sparked ideas and collaborations that continue to bear fruit in the mathematical sciences.

At the Institut des Hautes \'Etudes Scientifiques (IHES) in Bures-sur-Yvette, Pierre Cartier was more than just a researcher; he was a pillar of the institution, significantly contributing to its prominence. For those who had the privilege of knowing him, he was more than a colleague--he was a mentor, a guide, always ready to illuminate the darkest paths of research, to answer the most complex questions, or to place a problem in its proper conceptual context.

His passion for sharing extended beyond theory. Pierre Cartier was an indefatigable traveler\footnote{Among his travels, he undertook five training missions to Vietnam, each lasting one to two months, in 1976, 1980, 1984, 1988, and 1997. He also visited Latin America, specifically Chile and Argentina, in 1986, 1988, 1995, 1996, and 1997. His journeys took him to the Czech Republic in January 1991, 1992, and 1993, as well as in September 1987. In September 1987, he also traveled to Romania, followed by Ukraine in September 1993. Between October and December 1993, he visited India, while Japan was part of his itinerary in April 1976 and again from October to December 1990. He spent time in Quebec during September and October 1992 in Montreal, and in March 1996, he visited Munich, Germany. Austria's Vienna was another destination in March 1993, and in May 1994, he traveled to St. Petersburg, Russia.}, putting his generous ideas into practice by teaching around the world, particularly in Vietnam, where he frequently visited to impart his knowledge. He also loved political discussions, particularly with Laurent Lafforgue, and although their opinions often differed, they found mutual pleasure in these exchanges, reflecting Pierre's open-mindedness and intellectual vitality.

Indeed, Pierre maintained his open-mindedness and dedication to engaging with and sharing ideas within the mathematical community throughout his life. In January 2020, Pierre delivered his last lecture, which took place at IHES. By that time, he preferred not to give lectures while standing for extended periods, so we adapted the format to make it more comfortable for him. 
To support this, we prepared slides based on a video of a similar lecture he had given in 2018. The idea was to play the video alongside the slides, allowing Pierre to pause it at any point to provide additional insights and explanations. The topic was fascinating, exploring a paper by Hiroshi Umemura and an idea by Yuri Manin on deriving quantum groups from Galois theory, framed in such a way that the Galois group would naturally manifest as a quantum group. The lecture went beautifully. Pierre answered questions, engaged in discussions, and paused as needed to further elaborate on the material. This final lecture in January 2020 was a testament to his enduring passion and enthusiasm, and he continued to enjoy being a part of the IHES community.

\vspace{0.3cm}

As we reflect on Pierre Cartier's extraordinary career, we celebrate not only his mathematical achievements but also the spirit with which he approached his work. His contributions have not only advanced the frontiers of knowledge but have also embodied the essence of what it means to be a mathematician: curious, rigorous, creative, and deeply connected to the broader intellectual landscape.

His life was entirely dedicated to the pursuit of knowledge, the exploration of mathematical mysteries, and service to the scientific community. His influence was considerable, and his legacy will continue to live on through his contributions, writings, and the inspiration he imparted to so many of us. Even in the face of illness, he displayed exemplary courage, never allowing his physical condition to diminish his intellectual vitality.

His example encourages us to pursue research and the dissemination of knowledge with the same passion and dedication that he embodied throughout his life. As Jean Cocteau wrote, ``The true tomb of the dead is the heart of the living." Pierre will always remain with us, in our hearts and minds, forever inspiring us with his immense contribution to the world of mathematics.

\vspace{0.3cm}

Pierre Cartier had a great fondness for metaphors. It is interesting to note that his quotes, describing his journey and love for traveling to teach mathematics ``elsewhere"--not necessarily his own mathematics, but mathematics in general--serve as metaphors for his way of life and his joy in sharing his passion for the subject:

\begin{quote} I could describe myself as a mathematician without borders, borrowing from a well-known saying. By this, I mean crossing boundaries, which allowed me to do mathematics in some rather remarkable countries--Vietnam, Iraq, Kurdistan, and others. Teaching mathematics in such places made the effort worthwhile.

Why is it interesting to cross borders? Because on the other side, things are different. It's always exciting to venture to the other side of the fence, to see what lies in the shade. What may seem uninteresting on one side can be a treasure on the other, offering a fresh perspective. Something that might seem trivial here could be significant there.

So yes, I know a lot about crossing boundaries! They are meant to be crossed! From a scientific perspective, the start of my career took me across many frontiers--I began as a radio astronomer and, after a few detours into philosophy, ended up as a mathematician.

To do good science, you need that: constant imagination. No prejudice, and as I've learned from experience, no fear that your ideas might be silly\footnote{Personal communication}. 
\end{quote}

\section{Acknowldgement}

We want to warmly thank Mark C. Wilson for his support and suggestions during the preparation of the tribute. We would like to acknowledge the refree's careful reading  and for their suggestion regarding Cartier's theorem.


\begin{thebibliography}{99}

\bibitem{notices}
Cartier P., {\it Notice sur les travaux scientifiques de Pierre Cartier des ann\'ees 1970 -- 1990}.

\bibitem{cartier3}
Cartier P. (1958), {\it Dualit\'e des vari\'et\'es ab\'eliennes}, (s\'eminaire bourbaki, 1958), par Pierre Cartier. \url{http://www.numdam.org/item/SB_1956-1958__4__379_0.pdf}

\bibitem{dixmier}
Cartier P., Dixmier J. (1958), {\it Vecteurs analytiques dans les repr\'esentations de groupes de Lie}. Amer. J. Math., 80 (1958), 131-145.

\bibitem{cartier1}
Cartier P. (1957), {\it A new operation on differential forms par Pierre Cartier}, note pr\'esent\'ee par J. Hadamard, Comptes rendus de l'Acad. des Sc. {\it Une nouvelle op\'eration sur les formes diff\'erentielles}, C. R. Acad. Sci. Paris, 244 (1957), 426-428. \href{https://gallica.bnf.fr/ark:/12148/bpt6k31965/f426.item}{https://gallica.bnf.fr/ark:/12148/bpt6k31965/f426.item} and \href{https://www.math.stonybrook.edu/~jiahao/Notes/cartier.pdf}{https://www.math.stonybrook.edu/\~\,jiahao/Notes/cartier.pdf}.

\bibitem{cartier2}
Cartier P. (1958), {\it Questions de rationalit\'e des diviseurs en g\'eom\'etrie alg\'ebrique}.  \url{http://www.numdam.org/item/10.24033/bsmf.1503.pdf}


\bibitem{cartier4}
Cartier P. (1966) {\it Fonctions polylogarithmes, nombres polyz\^{e}tas et groupes pro-unipotents}, S\'eminaire Bourbaki, no. 885.

\bibitem{cartier5}
Cartier P. (1969), {\it Groupes formels, fonctions automorphes et repr\'esentations galoisiennes}, S\'eminaire Bourbaki, vol. 1968/69, exp. no. 347, 97--127.

\bibitem{cartier6}
Cartier P. (1971), {\it Probl\`emes math\'ematiques de la th\'eorie quantique des champs}, S\'eminaire N. Bourbaki, 1971, exp. no 388, p. 107-122

\bibitem{cartier7}
Cartier P. (1972), {\it Groupes alg\'ebriques et groupes formels.} S\'eminaire de l'Institut des Hautes \'etudes Scientifiques.


\bibitem{cartier8}
Cartier P. (1974), {\it Probl\`emes math\'ematiques de la th\'eorie quantique des champs II : prolongement analytique}, S\'eminaire N. Bourbaki, 1974, exp. no 418, p. 1-33


\bibitem{cartier10}
Cartier P. (2003), {\it A mad day's work: from Grothendieck to Connes and Kontsevich. The evolution of concepts of space and symmetry}, Bulletin of the American Mathematical Society, 38(4), 389-408.

\bibitem{connes}
Connes A., Kreimer D. (1998), {\it Hopf algebras, renormalization and non-commutative geometry}, Comm. Math. Phys., 199, 203-242.

\bibitem{cartier11}
Cartier P., DeWitt-Morette C. (2004), {\it Functional Integration : Action and Symmetries}. (Cambridge University Press, 2004). 

\bibitem{cartier12}
Cartier P. (2006), {\it Introduction aux alg\`ebres de Hopf}, sept. 2006, IH\'ES.  \url{https://people.math.osu.edu/kerler.2/VIGRE/InvResPres-Sp07/cartier-IHES.pdf}

\bibitem{cartier010}

Cartier P., Patras F. (2021), {\it Classical Hopf algebras and their application}, Springer.

\bibitem{cartier109}
Cartier P. (1962), {\it Groupes alg\'ebriques et groupes formels}, pp. 87--111 in : {\it Colloque sur la th\'eorie des groupes alg\'ebriques, (Bruxelles, 1962)}, Gauthier-Villars, 1962, p.109.


\bibitem{cartier01}
Cartier, P. (1963), {\it A Generalization of the Notion of Distribution}. In S\'eminaire Schwartz 1962/63 Expos\'e No. 7. Paris: \'Ecole Normale Sup\'erieure. 
    
  \bibitem{cartier02}
  Cartier, P. (1966), {\it Une Nouvelle D\'efinition des Distributions}. In S\'eminaire Schwartz 1965/66 Expos\'e No. 17. Paris: \'Ecole Normale Sup\'erieure.
    
\bibitem{cartier03}
Cartier, P., Schwartz, L. (1967), {\it Th\'eorie des Distributions et D\'eveloppement Asymptotique}. In Annales de l'Institut Fourier, 17(1), 293-364.
    
\bibitem{cartier04}
Cartier, P. (1967), {\it Fonctions g\'en\'eratrices, distributions et transformations de Fourier}. In S\'eminaire de Th\'eorie des Nombres (1965/1966), 109-141. Paris: \'Ecole Normale Sup\'erieure.
    
\bibitem{cartier05}
Cartier, P. (1968), {\it Integration Over Finite and Infinite-Dimensional Spaces}. In Proceedings of the International Congress of Mathematicians (Moscow), 421-436.
    
\bibitem{cartier06}
Cartier, P., Schwartz, L. (1971), {\it Th\'eorie des Distributions. In S\'eminaire Bourbaki}, 23e ann\'ee (1970/1971), Exp. No. 399. Paris: Soci\'et\'e Math\'ematique de France.
    
\bibitem{cartier07}
Cartier, P. (1980), {\it On the Structure of Distributions and Their Applications. In Lecture Notes in Mathematics}, Vol. 755, pp. 38-55. Berlin, Heidelberg: Springer.
    
\bibitem{cartier08}
Cartier, P. (1983), {\it Distribution Theory and Probability Theory: An Analytical Approach}, In S\'eminaire de Probabilit\'es XVII 1981/82, Lecture Notes in Mathematics, Vol. 986, pp. 23-47. Berlin, Heidelberg: Springer.

\bibitem{CK} A. Connes, D. Kreimer, {\it Renormalization in quantum
   field theory and the Riemann-Hilbert problem. I. The Hopf algebra
   structure of graphs and the main theorem}. Comm. Math. Phys. 210 (2000),
   no. 1, 249--273.

 \bibitem{CM} A. Connes, M. Marcolli (2008),  { \it Noncommutative Geometry, Quantum Fields, and Motives}, Colloquium Publications, Vol.55, American Mathematical Society.


\bibitem{lang}
Lang S. (1962), {\it Transcendental points on group varieties, Topology}, 1 (1962), p. 313-318.


\bibitem{fresan}
Fresan J. (2009), {\it The Castle of Groups}, \href{http://javier.fresan.perso.math.cnrs.fr/castle-of-groups.pdf}{http://javier.fresan.perso.math.cnrs.fr/castle-of-groups.pdf}.


\bibitem{cartier13}
Cartier P. (2007), {\it L'universalisme math\'ematique}, Diog\`ene 3 no 219, pp. 82-94 

\bibitem{cartier15}
Cartier P.,  Connes A., Lafforgue L. \& al. (2021), {\it Lectures grothendieckiennes}, eds. F. Ja\"{e}ck, SMF, Hors collection.
\end{thebibliography}
\end{document}